\documentclass[10pt]{article}
\usepackage{amscd,latexsym,amsmath,amssymb,amsfonts,graphicx,fancyhdr,setspace,xypic}
\usepackage{latexsym,amsmath,amssymb,amsfonts,verbatim}
\usepackage{hyperref}
\usepackage{amsmath,amscd}
\usepackage{chngcntr}

\author{
Junwu Tu\thanks{Mathematics Department, University of Oregon, Eugene OR 97403, USA, {\em e-mail:}{\tt junwut@uoregon.edu}}}
\title{ On the reconstruction problem in mirror symmetry}
\date{}

\counterwithout{paragraph}{subsection}
\counterwithin{paragraph}{section}

\DeclareFontFamily{U}{rsf}{}
\DeclareFontShape{U}{rsf}{m}{n}{
  <5> <6> rsfs5 <7> <8> <9> rsfs7 <10-> rsfs10}{}
\DeclareMathAlphabet{\mathscr}{U}{rsf}{m}{n}

\DeclareMathAlphabet{\mathgth}{U}{euf}{m}{n}

\DeclareFontFamily{U}{cyr}{}
\DeclareFontShape{U}{cyr}{m}{n}{
  <5> wncyr5 <6> wncyr6 <7> wncyr7 <8> wncyr8 <9> wncyr9 <10-> wncyr10}{}
\DeclareMathAlphabet{\mathcyr}{U}{cyr}{m}{n}

\input cyracc.def

\makeatletter
\def\operator@font{\sf}
\makeatother

\newcommand{\gog}{\mathgth{g}}

\newcommand{\cF}{{\mathscr F}}

\newcommand{\cI}{{\mathscr I}}
\newcommand{\cJ}{{\mathscr J}}
\newcommand{\cM}{{\mathscr M}}

\newcommand{\cO}{{\mathscr O}}

\newcommand{\cU}{{\mathscr U}}

\newcommand{\can}{{\mathsf{can}}}
\newcommand{\dR}{{\mathsf{dR}}}

\newcommand{\MC}{\mathsf{MC}}

\newcommand{\chk}{{\scriptscriptstyle\vee}}
\newcommand{\R}{\mathbb{R}}

\newcommand{\val}{{\mathsf{val}}}

\DeclareMathOperator{\Spec}{Spec}

\DeclareMathOperator{\Hom}{Hom}

\DeclareMathOperator{\GL}{GL}

\DeclareMathOperator{\id}{id}
\DeclareMathOperator{\ev}{ev}

\newcommand{\T}{\mathbb{T}}
\newcommand{\B}{\mathbb{B}}
\newcommand{\Coder}{\mathsf{Coder}}

\newcommand{\inter}{\mathsf{int}}
\newcommand{\ext}{\mathsf{ext}}

\newcommand{\ra}{\rightarrow}

\newcommand{\C}{\mathbb{C}}

\newcommand{\Z}{\mathbb{Z}}

\newcommand{\remark}{\noindent\textbf{Remark: }}
\newcommand{\proof}{\noindent\textbf{Proof. \;}}

\newtheorem{theorem}{Theorem}[section]

\newtheorem{lemma}[theorem]{Lemma}
\newtheorem{assumption}[theorem]{Assumption}
\newtheorem{corollary}[theorem]{Corollary}

\newtheorem{definition}[theorem]{Definition}
\newtheorem{definition-theorem}[theorem]{Definition-Theorem}

\renewcommand{\phi}{\varphi}

\makeatletter
\newcommand{\qed}{\hfill \ensuremath{\Box}}
\newcommand*{\doublerightarrow}[2]{\mathrel{
  \settowidth{\@tempdima}{$\scriptstyle#1$}
  \settowidth{\@tempdimb}{$\scriptstyle#2$}
  \ifdim\@tempdimb>\@tempdima \@tempdima=\@tempdimb\fi
  \mathop{\vcenter{
    \offinterlineskip\ialign{\hbox to\dimexpr\@tempdima+1em{##}\cr
    \rightarrowfill\cr\noalign{\kern.01ex}
    \rightarrowfill\cr}}}\limits^{\!#1}_{\!#2}}}
\newcommand*{\triplerightarrow}[1]{\mathrel{
  \settowidth{\@tempdima}{$\scriptstyle#1$}
  \mathop{\vcenter{
    \offinterlineskip\ialign{\hbox to\dimexpr\@tempdima+1em{##}\cr
    \rightarrowfill\cr\noalign{\kern.01ex}
    \rightarrowfill\cr\noalign{\kern.01ex}
    \rightarrowfill\cr}}}\limits^{\!#1}}}
\newcommand{\rnum}[1]{\romannumeral #1}
\newcommand{\Rnum}[1]{\expandafter\@slowromancap\romannumeral #1@}
\makeatother
\numberwithin{equation}{section}
\begin{document}
\maketitle
\begin{abstract}

Let $\pi:M\ra B$ be a Lagrangian torus fibration with singularities such that the fibers are of Maslov index zero, and unobstructed. The paper constructs a rigid analytic space $M^\chk_0$ over the Novikov field which is a deformation of the semi-flat complex structure of the dual torus fibration over the smooth locus $B_0\subset B$ of $\pi$. Transition functions of $M^\chk_0$ are obtained via $A_\infty$ homomorphisms which captures the wall-crossing phenomenon of moduli spaces of holomorphic disks.

\end{abstract}
\setcounter{secnumdepth}{4}
\section{Introduction}

\paragraph{Backgrounds.} Let $(M,\omega)$ be a symplectic manifold endowed with a (special) Lagrangian torus fibration $\pi: M\ra B$ possibly with singularities.  Associated to the fibration $\pi$ (plus additional polarization data), the reconstruction problem in mirror symmetry concerns about constructing a complex manifold $M^\chk$ which conjecturally should be a mirror partner for $M$. As such this problem lies at the heart towards a mathematical understanding of the mysterious mirror duality.

The reconstruction problem was first discussed by Kontsevich and Soibelman in~\cite[Section 7.1 Remark 19]{KS}. Later on in a series of remarkable papers~\cite{GS},~\cite{GS2},~\cite{GS3} Mark Gross and Bernd Siebert initiated, and studied in depth this problem from a more combinatorial point of view (a global version of toric geometry).
There are also nice concrete examples where this reconstruction process is explicitly realized, see for instance~\cite{Auroux} by Auroux, and~\cite{CLL} by Chan-Lau-Leung.

In the mean time Fukaya, Oh, Ono, and Ohta are developing Lagrangian Floer theory and its obstruction theory~\cite{FOOO}. After their work it has become evident that the homotopy theory of $A_\infty$ algebras is the right language for Lagrangian Floer theory. For example, associated to a (spin) Lagrangian submanifold $L$ in a symplectic manifold $M$ one can define an $A_\infty$ algebra structure on the de Rham complex of $L$ with coefficients in certain Novikov ring. This construction involves various choices that makes this $A_\infty$ structure not an invariant of symplectic geometry of the pair $(M,L)$. However the upshot is that the homotopy type of this structure is a symplectic invariant! 

Later on in~\cite{Fukaya}, inspired by Lagrangian Floer theory, Fukaya  briefly outlined a more conceptual approach to understand the reconstruction problem via gluing of Maurer-Cartan moduli spaces associated $A_\infty$ algebras. It is expected that once this construction is properly understood, homological mirror symmetry conjecture will follow from it with less efforts (see~\cite{Tu} for the local situation).

\paragraph{The current work.} In this paper we carry out Fukaya's approach to deal with the reconstruction problem over the smooth locus $B_0\subset B$ of a Lagrangian torus fibration $\pi: M\ra B$ under certain assumptions. Our main result is the following theorem. We refer to Corollary~\ref{coro:main} in Section~\ref{sec:re} for a more precise formulation of it.

\begin{theorem}
\label{intro:main}
Assume that Lagrangian torus fibers over $B_0$ are unobstructed, and are of Maslov index zero in $(M,\omega)$. Then there is a natural construction of a rigid analytic space $M^\chk_0$ over the Novikov field $\Lambda$ fibered over $B_0$, which gives a deformation of Kontsevich-Soibelman's construction~\cite[Section 7, Definition 22]{KS}.
\end{theorem}
\medskip
\remark The idea to use non-Archimedean geometry to deal with possible convergence issues in Floer theory is due to Kontsevich and Soibelman~\cite{KS}. Namely the mirror manifold $M^\chk$ to be constructed should not be a complex manifold over $\C$, but rather a rigid analytic space over a valuation ring. 

\medskip
\remark The explicit gluing formulas~\ref{eq:Psi} in the construction of $M^\chk_0$ involves certain instanton corrections from symplectic geometry. If the instanton corrections were not presented, then these formulas reduces to the one used in~\cite{KS}. We also note that for the construction of $M^\chk_0$ we do not need a polarization data on $M$. Such a polarization seems to be related to constructing a compatification of the space $M^\chk_0$. At present we do not know how to do this.

\medskip
\remark The Maslov index zero condition is automatic for \emph{special} Lagrangian torus fibrations in a Calabi-Yau manifold. However the unobstructedness assumption is not known in that case, even though expected. On the other hand, in view of~\cite{Tu}, the unobstructedness assumption is necessary for the purpose of homological mirror symmetry. Indeed without this assumption the Maurer-Cartan moduli spaces involved carry too less information in order to have an equivalence predicted by homological mirror symmetry.

\paragraph{Acknowledgement.} I am grateful to Kenji Fukaya for his encouragement on the subject as well as sharing a slide of his talk at MSRI. I also thank Vadim Vologodsky for useful discussions on non-Archimedean geometry. Finally thanks to the mathematics department of University of Oregon for providing excellent research condition.

\newpage

\section{Algebraic framework}
\label{sec:algebra}

In this section we recall certain homotopy theory of gapped and filtered $A_\infty$ algebras developped by Fukaya in~\cite{Fukaya}. Another reference is Fukaya's joint work with Oh, Ohta and Ono~\cite{FOOO}. The construction of canonical models in the gapped and filtered context was first discussed in~\cite{FOOO2}, see also~\cite{FOOO} and~\cite{Fukaya}.

\paragraph{Monoids.} Since our primary geometric applications concern only Lagrangians with vanishing Maslov index, we shall work with a monoid $G\subset \R^{\geq 0}$ which keeps track of the energy of pseudo-holomorphic curves~\footnote{In~\cite{Fukaya} and~\cite{FOOO} the authors used monoids in $\R^{\geq 0}\times 2\Z$ in which the component $2\Z$ keeps track of the Maslov index.}. In view of Gromov's compactness theorem we require that 
\[ |G\cap [0,E]| < \infty, \;\;\forall E\in \R^{\geq 0}.\] 
Given such a monoid $G$ we can form a ring $\Lambda_G^\C$ consisting of formal sums of the form $\sum_{i=0}^\infty a_i T^{\beta_i}$ such that $ a_i\in \C, \beta_i\in G$. The finiteness condition above implies that 
\[|\left\{ i\mid a_i\neq 0, \beta_i\leq E\right\}| < \infty, \;\;\forall E\in \R^{\geq 0}.\]
We define a valuation on $\Lambda_G^\C$ by 
\[ \val (\sum_{i=0}^\infty a_i T^{\beta_i}):= \inf_{\left\{ i\mid a_i\neq 0\right\}} \beta_i.\] 
This valuation map induces a filtration on $\Lambda_G^\C$ by setting
\[ F^E(\Lambda_G^\C):=\left\{ a\in \Lambda_G^\C\mid \val(a)\geq E\right\}\]
for a fixed energy $E\in \R^{\geq 0}$. This filtration defines a topology on the ring $\Lambda_G^\C$, making it a topological ring (i.e. its product is continuous with respect to the topology). Furthermore the ring $\Lambda_G^\C$ is complete with respect to this topology.

\paragraph{Gapped filtered $A_\infty$ algebras.} Let $(\overline{V},D,\circ)$ be a differential graded algebra over $\C$, and denote by $\B \overline{V}$ its cobar differential graded coalgebra. Explicitly $\B \overline{V}$ is the free graded coalgebra $\T (\overline{V}[1])$ endowed with a degree one coderivation $Q_0$ determined by two degree one linear maps
\begin{align*}
D&: \overline{V}[1]\ra \overline{V}[1] \\
\circ &: \overline{V}[1]\otimes\overline{V}[1]\ra \overline{V}[1].
\end{align*}
Consider a monoid $G$ as explained above, and form the tensor product
\[ V:=\overline{V}\hat{\otimes} \Lambda_G^\C.\]
Here $\hat{\otimes}$ stands for first taking ordinary tensor product, and then taking topological completion with respect to the topology induced from the filtration mentioned above. Explicitly elements of $V$ are of the form $\sum_{i=0}^\infty v_i T^{\beta_i}$ where $v_i\in \overline{V}$ and $\beta_i\in G$. These series are required to satisfy the same finiteness condition as in the definition of $\Lambda_G^\C$.

The linear space $V$ also has a valuation map defined in the same way as that of $\Lambda_G^\C$. With the induced topology the $\Lambda_G^\C$-module structure on $V$ is continuous. 

For later usage we shall call elements with strictly positive valuations positive. For example, positive elements of $V$ are formal series such that the coefficient of $T^0$ is zero.

We extend the differential $D$ and the product $\circ$ on $\overline{V}$ to the whole space $V$ by $\Lambda_G^\C$-linearity to obtain a differential graded algebra structure on $V$ over $\Lambda_G^\C$. Denote by $\B V$ its cobar construction, and again by $Q_0$ its cobar differential. 

\begin{definition}
\label{def:a-infinity}
A $G$-gapped filtered $A_\infty$ algebra structure on $V$ is a \textbf{degree one, positive}  $\Lambda_G^\C$-linear coderivation $\delta$ on the coalgebra $\B V:=\T(V[1])$ such that 
\[(Q_0+\delta)^2 =0\]
\end{definition}
\remark Euivalently the element $\delta$ maybe considered a \emph{positive} Maurer-Cartan element of the differential graded Lie algebra $\Coder(\B V)$ where we consider the standard Lie structure on $\Coder(\B V)$, and its differential is given by $[Q_0,-]$. Indeed sine the differential $Q_0$ on $\Coder(\B V)$ squares to zero, the equation $(Q_0+\delta)^2=0$ is equivalent to $[Q_0,\delta]+\frac{1}{2}[\delta,\delta]=0$ which is precisely the Maurer-Cartan equation of $\Coder (\B V)$.

\paragraph{Explicit formula for $A_\infty$ algebras.} Let us unwind the definition of a $G$-gapped filtered $A_\infty$ structure in terms of multi-linear maps on $V$. For this note that there is a bijection between the two sets 
\[\Coder(\B V)\leftrightarrow \Hom(\B V, V)=\prod_{k=0}^\infty \Hom(V[1]^k, V[1])\] 
which sends an element $\phi\in \Coder(\B V)$ to $\pi \circ\phi$ where $\pi: \B V\ra V[1]$  is the canonical projection map. The backward map is given by extending a given multi-linear map in $\Hom(V[1]^k, V[1])$ by the co-Leibniz rule to a coderivation on $\B V$. Denote by $m_k$ the component in $\Hom(V[1]^k, V[1])$ corresponding to the coderivation $(Q_0+\delta)$ in Definition~\ref{def:a-infinity}. The continous $\Lambda^\C_G$-linear map $m_k: V[1]^k\ra V[1]$ can written as
\[ m_k=\sum_{\beta\in G} m_{k,\beta} T^\beta \]
for $\C$-linear maps $m_{k,\beta}: \overline{V}[1]^k\ra \overline{V}[1]$ of degree one. The positivity of $\delta$ implies that 
\begin{equation}
\label{eq:zero}
m_{k,0}=\begin{cases}
D & k=1 \\
\circ & k=2 \\
0 & \mbox{otherwise.}
\end{cases}
\end{equation}
Written in terms of $m_{k,\beta}$ the equation $(Q_0+\delta)^2$ gets translated into the equations
\begin{equation}
\label{eq:ainfinity}
 \sum_{i+j+k=N, i\geq 0, j\geq 0, k\geq 0}\sum_{\beta_1+\beta_2=\beta} m_{i+k+1,\beta_1}(\id^i\otimes m_{j,\beta_2} \otimes \id^k) =0
\end{equation}
for all $N\geq 0, \beta\in G$. Here Koszul sign convention is assumed when forming tensor products of linear maps.

\paragraph{Pseudo-isotopies of $A_\infty$ algebras.} Let $\Omega^*_{\Delta^1}$ be the set of $C^\infty$ differential forms on the $1$-simplex $\Delta^1=[0,1]$ that are constant in a neighborhood of $0$ and $1$. The boundary condition is due to a technical reason: being able to perform gluing of pseudo-isotopies in the smooth category.
\begin{definition}
A pseudo-isotopy between two $G$-gapped filtered $A_\infty$ algebra structures $\delta_0$ and $\delta_1$ on $V$ is given by a positive Maurer-Cartan element $\gamma\in \MC(\Coder(\B V) \hat{\otimes} \Omega^*_{\Delta^1})$. The element $\gamma$ is required to satisfy the boundary conditions
\[ i_0^*\gamma=\delta_0, \mbox{\;\;\; and\;\;\;\;} i_1^*\gamma=\delta_1.\]
Here the maps $i_j^*$ denote pull-backs of differential forms via the inclusions $i_j: j\ra [0,1]$ for $j=0, \mbox{ or }1$.
\end{definition}

\begin{lemma}
\label{lem:equiv}
Pseudo-isotopy defines an equivalence relation on the set of all $G$-gapped filtered $A_\infty$ structures on $V$.
\end{lemma}

\proof Let $\delta$ be an $A_\infty$ algebra structure, then $\gamma=\delta$ viewed as a ``constant" element of $\Coder(\B V)\otimes \Omega^*_{\Delta^1})$ is a pseudo-isotopy between $\delta$ and itself. Reversing a given pseudo-isotopy $\gamma^t\mapsto \gamma^{1-t}$ proves reflexive property. For the transitivity let $\gamma_1$ be a pseudo-isotopy between $\delta_0$ and $\delta_1$. Let $\gamma_2$ be another pseudo-isotopy between $\delta_1$ and $\delta_2$. There are the following two maps between intervals
\begin{align*}
\lambda_1: [0,1/2] \ra [0,1], \;\;\; & t\mapsto 2t ;\\
\lambda_2: [1/2] \ra [0,1], \;\;\; & t\mapsto 2t-1.
\end{align*}
Then $\lambda_1^*\gamma_1$ and $\lambda_2^*\gamma_2$ are differential forms on $[0,1/2]$ and $[1/2,1]$ with values in the Lie algebra $\Coder(\B V)$. Moreover they are both $\delta_1$ in a neighborhood of $1/2$. Hence we can glue them to form a $C^\infty$ form on $[0,1]$ with values in $\Coder(\B V)$. We denote by the result of this gluing by $\gamma_1\sharp \gamma_2$ to mimic the ``concatenation of paths" in topology. Then $\gamma_1\sharp \gamma_2$ is a pseudo-isotopy between $\delta_0$ and $\delta_2$. The lemma is proved.\qed


\paragraph{Explicit formula for pseudo-isotopies.} Let $\delta_0$, $\delta_1$ be two $A_\infty$ structures on $V$, and denote by $m^0_k$, $m^1_k$ the corresponding multi-linear maps on $V[1]$. A pseudo-isotopy $\gamma$ between $\delta_0$ and $\delta_1$ is a Maurer-Cartan element of the form 
\[ \gamma=\delta^t+ h^tdt\]
for some positive elements $\delta^t, h^t \in \Coder(\B V)\hat{\otimes} C^\infty([0,1])$ such that $\delta^t$ is of degree one, $h^t$ is of degree zero, satisfying the initial condition 
\[ \delta^t|_{t=0}=\delta_0, \mbox{\;\; and \;\;\;} \delta^t|_{t=1}=\delta^1.\] 
Moreover the Maurer-Cartan equation for $\gamma$ implies that
\[ \left\{
  \begin{array}{l}
     [Q_0+\delta^t,Q_0+\delta^t]=0,\\
     \frac{d\delta^t}{dt}=[Q_0+\delta^t,h^t].\\
  \end{array} \right. \]
Let $m^t_{k,\beta}$ and $h^t_{k,\beta}$ be the associated multi-linear maps on $\overline{V}[1]$ corresponding to the coderivations $Q_0+\delta^t$ and $h^t$. The positivity of $\gamma=\delta^t+h^tdt$ implies that $m^t_{k,0}$ is as in equation~\ref{eq:zero} for all $t\in [0,1]$, and $h^t_{k,0}\equiv 0$. In terms of these component maps, the first equation above says that $m^t_k=\sum_\beta m^t_{k,\beta}$ form a family of $G$-gapped filtered $A_\infty$ structure on $V$. The second equation can be expanded to get
\begin{multline}
\label{eq:isotopy}
\frac{dm^t_{N,\beta}}{dt} =-\sum_{i+j+k=N}\sum_{\beta_1+\beta_2=\beta} m^t_{i+k+1,\beta_1}(\id^i\otimes h^t_{j,\beta_2} \otimes \id^k)
\\+\sum_{i+j+k=N} \sum_{\beta_1+\beta_2=\beta} h^t_{i+k+1,\beta_1}(\id^i\otimes m^t_{j,\beta_2} \otimes \id^k) 
\end{multline}
where the equation holds for all $N\geq 0$, and  $\beta \in G$.
\paragraph{From pseudo-isotopies to $A_\infty$ homomorphisms.} Recall a filtered $A_\infty$ homomorphism from $(A,\delta_A)$ to $(B,\delta_B)$ is a continuous map of differential graded coalgebras $F: \B A \ra \B B$. If $m^A_{k,\beta}$ and $m^B_{k,\beta}$ are the corresponding structure constants for $\delta_A$ and $\delta_B$, the map $F$ can be realized as multi-linear maps $F_{k,\beta}:(\overline{A}[1])^k\ra \overline{B}[1]$ of degree zero satisfying 
\begin{multline}
\label{eq:homo}
\sum_{0\leq j,i_1+\cdots+i_j=N}\sum_{\beta_0+\beta_{i_1}+\cdots+\beta_{i_j}=\beta}m_{j,\beta_0}^B(F_{i_1,\beta_{i_1}}\otimes\cdots\otimes F_{i_j,\beta_{i_j}}) =\\ =\sum_{i+j+k=N}\sum_{\beta_1+\beta_2} F_{i+k+1,\beta_1}(\id^i\otimes m_{j,\beta_2}^A\otimes \id^k)
\end{multline}
where the equation holds for all $N\geq 0$, and $\beta\in G$. 

Given a pseudo-isotopy $\gamma$ between $(V,\delta_0)$ and $(V,\delta_1)$ we can construct an $A_\infty$ homomorphism from $(V,\delta_0)$ to $(V,\delta_1)$, which we now explain. We need a few ingredients for this construction.

\paragraph{Stable trees with decorations.} Let $k\geq 0$ be a nonnegative integer. Consider a rooted ribbon tree $T$ with $k$ leaves. Here the word ``ribbon" means that, for every vertex of the tree, a cyclic ordering is assigned to the set of half edges incident to that vertex. The root of a tree $T$ will be denoted by $r_T$, the set of edges by $E(T)$, and the set of vertices by $V(T)$. We remark that the set $V(T)$ is the disjoint union of $V^\inter(T)$ consisting of interior vertices and $V^\ext(T)$ consisting of $k+1$ exterior vertices (corresponding to the root and the $k$ leaves). The exterior vertices are necessarily of valency $1$, for the interior vertices we allow them to have arbitrary positive valency. A $G$-decoration on $T$ is a map $\lambda: V^\inter(T) \ra G$. We call a $G$-decorated ribbon tree stable if every interior vertex decorated by $0\in G$ has at least valency $3$. We denote by $\cO(k,\beta)$ be the set of stable $G$-decorated rooted ribbon trees with $k$ leaves such that
\[ \sum_{v\in V^\inter(T)} \lambda(v)=\beta.\]
\begin{lemma}
\label{lem:finite}
The set $\cO(k,\beta)$ is a finite set for all $k\geq 0$, and $\beta\in G$.
\end{lemma}

\proof It suffices to show that the number of the set $V^\inter(T)$ is finite. By the finiteness assumption on the monoid $G$ this is equivalent to show that the set of interior vertices decorated by $0\in G$ is finite. For this we observe that since $T$ is a tree which is contractible, hence its Euler characteristic is
\[ (k+1)+|V^\inter(T)|- |E(T)|=1.\]
Denote by $x$ the number of vertices decorated by zero, and by $y$ the rest interior vertices. So we have $x+y=|V^\inter(T)|$. By the stableness assumption there are at least $3x$ half edges incident to these vertices. Thus $|E(T)|\geq 3x/2$, together with the above qualities implies that $x\leq 2(k+y)$. Since the right hand side is bounded, $x$ is also bounded. The lemma is proved.
\qed

\medskip
Let $(T,\lambda) \in \cO(k,\beta)$ be a stable $G$-decorated tree. We define a partial ordering on the set $V^\inter(T)$ as follows. We say $v_1\leq v_2$ if $v_2$ is contained in the shortest path from $v_1$ to the root $r_T$. Here we measure a path by its number of edges in $T$. Let $t\in[0,1]$, a time ordering on $T$ bounded by $t$ is a non-decreasing map $\tau: V^\inter(T) \ra [0,t]$. Since this partial ordering only depends on the tree $T$, and not on the decoration $\lambda$, we denote by $\cM^t(T)$ the set of all time orderings on $T$. We endow it with the subset topology of $[0,t]^{|V^\inter(T)|}$. As such it is compact. The Euclidean measure on $[0,t]^{|V^\inter(T)|}$ induces a measure on $\cM^t(T)$ which we denote by $d\tau$.

For each triple $(T,\lambda,\tau)$ we will define a multi-linear map $\rho(T,\lambda,\tau):\overline{V}[1]^k\ra \overline{V}[1]$ of degree zero. On each internal vertex $v\in V^\inter(T)$ we put the linear map $h^{\tau(v)}_{\val(v)-1,\lambda(v)}$ where $\val(v)$ is the valency of the vertex $v$. Then the map $\rho(T,\lambda,\tau)$ is simply the ``operadic" composition of these linear maps according to the tree $T$. 

Finally we can define the structure maps  $F_{k,\beta}^t: \overline{V}[1]^k \ra \overline{V}[1]$ of an $A_\infty$ homomorphism by
\begin{equation}
\label{eq:phi}
F_{k,\beta}^t:= \sum_{(T,\lambda)\in \cO(k,\beta)} \rho^t(T,\lambda) \mbox{\;\;\;\; where\;\;\;} \rho^t(T,\lambda):=\int_{\cM^t(T)} \rho(T,\lambda,\tau) d\tau.
\end{equation}
By the previous lemma this is only a finite sum, hence is well-defined. The main property for these complicatedly defined maps is the following derivative formula.

\begin{lemma}
\label{lem:derivatives}
For each $(T,\lambda)\in \cO(k,\beta)$ we have the following formula
\[ \frac{d\rho^t(T,\lambda)}{dt}= h^t_{l,\lambda(v_0)} [\rho^t(T^{(1)},\lambda^{(1)})\otimes\cdots\otimes \rho^t(T^{(l)},\lambda^{(l)})].\] 
Here the vertex $v_0$ is the unique vertex connected the root of $T$, and the trees $T^{(1)},\cdots,T^{(l)}$ are obtained so that $\mu_l\circ(T^{(1)}\otimes\cdots\otimes T^{(l)})= T$ where $\mu_l$ denote the star graph with $l$ inputs. Namely if we glue the roots of $T^{(i)}$'s into the inputs of $\mu_l$, we obtain the tree $T$. Decorations $\lambda^{(1)},\cdots,\lambda^{(l)}$ are induced from that of $\lambda$.
\end{lemma}

\proof The integral $\int_{\cM^t(T)} \rho(T,\lambda,\tau) d\tau$, by its definition, may be rewritten as
\[ \int_0^t h^{\tau(v_0)}_{l,\lambda(v_0)}\circ[\rho^{\tau(v_0)}(T^{(1)},\lambda^{(1)})\otimes\cdots\otimes\rho^{\tau(v_0)}(T^{(l)},\lambda^{(l)})] d\tau(v_0).\]
The lemma then follows from fundamental theorem of calculus.\qed

\begin{lemma}
\label{lem:homo}
The maps $F_k^t$ defined by formula~\ref{eq:phi} form an $A_\infty$ homomorphism $F: (V,m^0)\ra (V,m^t)$.
\end{lemma}

\proof By formula~\ref{eq:homo} we need to compare the following two terms for all $n\geq 0$, and $\beta\in G$:
\begin{equation*}
\label{eq:I}
\Rnum{1}_{n,\beta}^t:= \sum_{i+k+j=n}\sum_{\beta_1+\beta_2=\beta} F_{i+j+1,\beta_1}^t(\id^i\otimes m_{k,\beta_2}^0 \otimes \id^j);
\end{equation*}
\begin{equation*}
\label{eq:II}
\Rnum{2}_{n,\beta}^t:=\sum_{i_1+\cdots+i_k=n}\sum_{\beta_0+\beta_{i_1}+\cdots+\beta_{i_k}=\beta} m^t_k(F_{i_1,\beta_{i_1}}^t\otimes\cdots\otimes F_{i_k,\beta_{i_k}}^t).
\end{equation*}
Denote by $\Delta_{n,\beta}^t:=\Rnum{1}_{n,\beta}^t-\Rnum{2}_{n,\beta}^t$ their difference. Computing $\frac{d\Delta_{n,\beta}^t}{dt}$ using formula~\ref{eq:isotopy} and Lemma~\ref{lem:derivatives} yields

\begin{align*}
 \frac{d\Delta_{n,\beta}^t}{dt}&= \sum_{i_1+\cdots+i_{j+l}+k=n} \sum_{\beta_0+\beta_1+\beta_{i_1}+\cdots+\beta_{i_{j+l}}=\beta}\\
& h^t_{j+l+1,\beta_0}(F_{i_1,\beta_{i_1}}^t\otimes\cdots\otimes F_{i_j,\beta_{i_j}}^t\otimes \Delta_{k,\beta_1}^t\otimes F_{i_{j+1},\beta_{i_{j+1}}}^t\otimes\cdots\otimes F_{i_{j+l},\beta_{i_{j+l}}}^t). \end{align*}
Moreover we have the initial conditions $\Delta^0_{n,\beta}=0$, and $\Delta^t_{n,0}=0$ for all $n\geq 0$, $\beta\in G$. 

In the case $n=0$, we can proceed the proof by doing induction on $\beta$. Indeed $\Delta^t_{n,0}$ serves as the initial step of the induction. For higher $\beta$ observe that the right hand side reduces to an ordinary differential equation by induction. Hence by uniqueness of solution of ordinary differential equations $\Delta_{0,\beta}^t\equiv 0$. In general we can proceed by induction also on the variable $n$, namely we assume $\Delta_{k,\beta}^t\equiv 0$ for all $k\leq n$ and $\beta<E$. Here $E$ is taken so that $E\in G$ and $(\beta,E)$ contains no other element in $G$. In this case we can show that $\Delta_{n,E}$ is zero which again follows from uniqueness of solution of ordinary differential equations. Thus the lemma is proved.\qed

\paragraph{Nerve of $\Coder(\B V)$.} Consider smooth differential forms $\Omega_{\Delta^n}^*$  on the standard $n$-simplex $\Delta^n:=\left\{(t_0,\cdots,t_n)\in \R^{n+1}| t_i\geq 0, \sum t_i =1.\right\}$ that are constant in a neighborhood of the vertices of $\Delta^n$. Let $V$ be as before, and denote by $\gog$ the graded Lie algebra $\Coder(\B V)$ consisting of continuous coderivations. We have considered the Maurer-Cartan elements of $\gog\otimes \Omega_{\Delta^0}^*=\gog$ which are $A_\infty$ structures on $V$, and that of $\gog \otimes \Omega_{\Delta^1}^*$ which are pseudo-isotopies between $A_\infty$ structures. 

The spaces $\MC(\gog)$ and $\MC(\gog\otimes\Omega_{\Delta^1}^*)$ are only ``pieces of an iceberg". Indeed the collection of differential graded Lie algebras $\left\{ \gog\otimes\Omega_{\Delta^n}^*\right\}_{n=0}^\infty$ form a simplicial differential graded Lie algebra $\gog_*$. Applying the Maurer-Cartan functor to it yields a simplicial set $\Sigma_*\gog:=\MC(\gog_*)$, which was introduced by Hinich~\cite{Hinich} and further investigated by Getzler~\cite{Getzler}.

Intuitively elements of $\Sigma_*\gog=\MC(\gog\otimes\Omega_{\Delta^n}^*)$ for $n\geq 2$ are pseudo-isotopies of pseudo-isotopies. We may think of elements of $\Sigma_0\gog$ as points, and elements of $\Sigma_1\gog$ as paths between these points.
\[ \gamma(0):=i_0^*\gamma=\delta_0\;\;\bullet \stackrel{\gamma}{\longrightarrow} \bullet\;\;\gamma(1):=i_1^*\gamma=\delta_1\]
Similarly elements of $\Sigma_2\gog$ are homotopies of paths, and elements of $\Sigma_3\gog$ are homotopies of homotopies of paths, and so on. In the following we shall call an element of $\Sigma_n \gog$ an $n$-pseudo-isotopy. For the purposes of this paper we only need to consider the cases for $n\leq 2$.

\paragraph{From $2$-pseudo-isotopies to homotopies of $A_\infty$ morphisms.} Let $\alpha\in \Sigma_2 \gog$ be a $2$-pseudo-isotopy, and let $\partial_0(\alpha)$, $\partial_1(\alpha)$, $\partial_2(\alpha)$ be its boundary pseudo-isotopies.
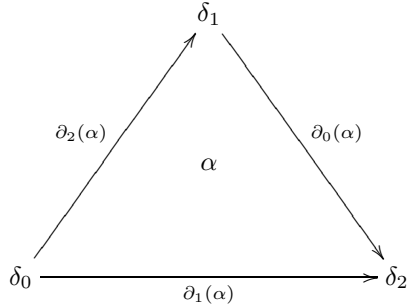
\[\begin{xy} 
(0,0)*+{\delta_0}="a"; (50,0)*+{\delta_2}="b";%
(25,35)*+{\delta_1}="c"; 
(25,15)*+{\alpha}="d";
{\ar@{->}^{\partial_2(\alpha)} "a";"c"};
{\ar@{->}_{\partial_1(\alpha)} "a";"b"}
{\ar@{->}^{\partial_0(\alpha)} "c";"b"}
\end{xy}\]
\begin{theorem}
\label{thm:homotopy}
The two $A_\infty$ morphisms $F[\partial_1(\alpha)]$ and $F[\partial_0(\alpha)]\circ F[\partial_2(\alpha)]$ obtained from boundary $1$-pseudo-isotopies of $\alpha$ are homotopic homomorphisms (see~\cite[Chapter 4]{FOOO} for the notion of homotopy used here).
\end{theorem}

\proof Let $u: [0,1]_s\times [0,1]_t \ra \Delta^2$ be an orientation preserving smooth map with the following boundary conditions:
\begin{align*}
u(0\times [0,1])&=0;\\
u(1\times [0,1])&=2;\\
u([0,1]\times 0)&= \partial_1;\\
u([0,1]\times 1)&=\partial_2 \cup \partial_0.
\end{align*}
Here on the right hand side of these equations the notations are illustrated in the figure above. Then the $t$-direction defines a homotopy between $F[\partial_1(\alpha)]$ and $F[\partial_0(\alpha)]\circ F[\partial_2(\alpha)]$. The lemmas is proved.\qed

\paragraph{Canonical models.} In the rest of this section we recall a few formulas for doing homological perturbation for $G$-gapped filtered $A_\infty$ algebras and their pseudo-isotopies. These formulas are important because we need to use their cyclic symmetry and exponential dependence in the next section.

Recall we have a differential graded algebra $(\overline{V},D,\circ)$ on the zero energy part of a $G$-gapped filtered $A_\infty$ algebra $V=\overline{V}\hat{\otimes} \Lambda^\C_G$. Let us denote by $\overline{H}$ the cohomolog of the complex $(V,D)$. Consider linear maps $i:\overline{H}\ra \overline{V}$, $p: \overline{V}\ra \overline{H}$ of degree zero, and $h: \overline{V}\ra \overline{V}$ of degree $-1$ satisfying conditions
\[ Di=0,\;\; pD=0,\;\; pi=\id,\;\; \id-ip=Dh+hD.\]
Such a triple $(i,p,h)$ is called a deformation retraction. Given a $G$-gapped filtered $A_\infty$ algebra $V$ and a deformation retraction $(i,p,h)$ on $\overline{V}$, we can construct a $G$-gapped filtered $A_\infty$ algebra structure on $H:=\overline{H}\otimes_\C \Lambda^\C_G$. This is known as the canonical model for $V$.

We describe this $G$-gapped filtered $A_\infty$ algebra on $H$ by providing explicit formulas for $m^\can_{k,\beta}$. For this let $(T,\lambda)\in \cO(k,\beta)$ be a $G$-decorated ribbon tree with $k$ leaves. We associate to each interior vertex $v\in V^\inter(T)$ the operator $m_{\val(v)-1,\lambda(v)}$, to each interior edge (edges that are not connected to a leave or the root) the homotopy operator $h$, to each leaves the embedding $i$, and finally to the unique edge incident to the root the projection map $p$. Define $\eta(T,\lambda)$ as the operadic composition of these multi-linear maps. Then the structure maps $m^\can_{k,\beta}$ on $H$ is defined as
\begin{equation}
\label{eq:can_1}
m^\can_{k,\beta}:= \sum_{(T,\lambda)\in \cO(k,\beta)} \eta(T,\lambda).
\end{equation}
Let $\gamma=\delta+hdt \in \Sigma_1\gog$ be a $1$-pseudo-isotopy on $V$, and let $m^t_{k,\beta}$, $h^t_{k,\beta}$ be its associated structure maps. We shall define $\gamma^\can$ which is a $1$-pseudo-isotopy on $H$. The formula for $m^{t,\can}_{k,\beta}$ is given by equation~\ref{eq:can_1} above allowing the parameter $t$ to enter. 

The formula for $h^{t,\can}_{k,\beta}$ is defined as follows. Consider $\cO^+(k,\beta)$ the set of triples $(T,\lambda,v)$ where $(T,\lambda)$ is a $G$-decorated ribbon trees with $k$ leaves in $\cO(k,\beta)$, and $v\in V^\inter(T)$ a chosen interior vertex of $T$. Then we associate operators to this triple $(T,\lambda,v)$ in the same way as in the construction of $m^\can_{k,\beta}$ except at the chosen vertex we put the operator $h^t_{\val(v)-1,\lambda(v)}$. Define $\mu(T,\lambda,v)$ to be the operadic composition, and the formula for $h^{t,\can}_{k,\beta}$
is defined by
\begin{equation}
\label{eq:can_2}
h^{t,\can}_{k,\beta}:= \sum_{(T,\lambda,v)\in \cO^+(k,\beta)} \mu(T,\lambda,v).
\end{equation}
\paragraph{Maurer-Cartan moduli spaces.} An important homotopy invariant associated to a $G$-gapped filtered $A_\infty$ algebra is its Maurer-Cartan moduli space. As before we use $m_{k,\beta}$, and $m_k:=\sum_\beta m_{k,\beta}$ for structure maps of a $G$-gapped filtered $A_\infty$ structure $\delta$ on $V$. Furthermore we denote by $\Lambda_G^+$ consisting of elements in $\Lambda_G^\C$ of positive valuation, and denote by $V^+:=\overline{V}\otimes\Lambda_G^+$.

Let $b\in V^+$ be a positive element of degree one. The Maurer-Cartan equation is defined as
\[ m_0+m_1(b)+m_2(b^2)+\cdots+m_k(b^k)+\cdots=0.\]
Note that the fact that $b$ has positive valuation ensures the convergence of this equation. Solutions of this equation are called Maurer-Cartan elements of $V$. Two Maurer-Cartan elements $b_0$ and $b_1$ are called gauge equivalent if there exists a Maurer-Cartan element $\theta\in V\otimes\Omega_{\Delta^1}^*$ such that $\theta$ is positive, and satisfies
\[ i_0^*\theta=b_0, \mbox{\;\; and \;\;} i_1^*\theta=b_1.\]
More explicitly the element $\theta$ is of the form $\theta=b(t)+c(t)dt$ for a one parameter degree one positive element $b(t)$, and one parameter degree zero element $c(t)$. They are subject to the following equations
\begin{equation}
\label{eq:mc}
\begin{cases}
m_0+m_1(b(t))+\cdots+m_k(b(t)^k)+\cdots=0, & \\
\frac{db(t)}{dt}+\sum_{k=0}^\infty\sum_{i+j=k}m_{k+1}(b(t)^i,c(t),b(t)^j)=0. &
\end{cases}
\end{equation}

\begin{lemma}
Gauge equivalence is an equivalence relation.
\end{lemma}

\medskip
\noindent
The lemma is proved in~\cite{FOOO}. The Maurer-Cartan moduli space $\MC^+(V,\delta)$~\footnote{Here we added a plus sign to indicate that we only consider elements of positive valuation.} of an $A_\infty$ algebra $(V,\delta)$ is then defined as the set of equivalence classes of Maurer-Cartan elements modulo gauge equivalence. The fundamental property of the assignment sending
\[ \mbox{$A_\infty$ algebra $(V,\delta)$} \mapsto \MC^+(V,\delta) \]
is that it factor through the homotopy category of $A_\infty$ structures on $V$. In particular homotopy equivalent $A_\infty$ structures have isomorphic Maurer-Cartan moduli spaces. More explicitly if $\Phi: A\ra B$ is a map of $G$-gapped filtered $A_\infty$ homomorphism, then the map
\begin{equation}
\label{eq:mc2}
 b\mapsto \Phi_0(1)+\Phi_1(b)+\Phi_2(b^2)+\cdots+\Phi_k(b^k)+\cdots \end{equation}
respects gauge equivalence, which descends to a map $\Phi_*: \MC^+(A)\ra \MC^+(B)$. In case when $\Phi$ is a homotopy equivalence this map is an isomorphism. In particular if $\gamma$ is a $1$-pseudo-isotopy between two $A_\infty$ algebra structures $\delta_0$ and $\delta_1$ on $V$, then the $A_\infty$ homomorphism $F(\gamma): (V,\delta_0) \ra (V,\delta_1)$ defined earlier induces an isomorphism on the associated Maurer-Cartan moduli spaces.

\section{Geometric realizations}
\label{sec:gr}

\paragraph{$A_\infty$ algebras associated to Lagrangian submanifolds.} Let $L$ be a relatively spin compact Lagrangian submanifold in a symplectic manifold $(M,\omega)$. In~\cite{FOOO} and~\cite{Fukaya} a curved $A_\infty$ algebra structure was constructed on $\Omega^*(L,\Lambda_0)$, the de Rham complex of $L$ with coefficients in certain Novikov ring $\Lambda_0$ over $\C$~\footnote{Strictly speaking this version of Floer theory was developed by K. Fukaya in~\cite{Fukaya} heavily based on the work of himself, Y.-G. Oh, H. Ohta and K. Ono~\cite{FOOO}.}. Here coefficient ring $\Lambda_0$ is defined by
\[ \Lambda_0:=\left\{\sum_{i=1}^\infty a_i T^{\lambda_i} \mid a_i\in \C, \lambda_i\in \R^{\geq 0}, \lim_{i\ra \infty} \lambda_i= \infty\right\}\]
where $T$ is a formal parameter of degree zero. To make connections with the ring $\Lambda^\C_G$ used in the previous section, we observe that there is an inclusion of monoids \[ G\hookrightarrow \R^{\geq 0},\]
which induces an inclusion of rings $\Lambda^\C_G\hookrightarrow \Lambda_0$. The $A_\infty$ algebra associated to $L$ can be constructed over $\Lambda^\C_G$ for some monoid $G$ depending on $L$. By extension of scalars we get an $A_\infty$ algebra over $\Lambda_0$. We shall refer to this $A_\infty$ algebra as the Fukaya ($A_\infty$) algebra of $L$.

Note that the map $\val:\Lambda_0\ra \R$ defined by $\val(\sum_{i=1}^\infty a_i T^{\lambda_i}):= \inf_{a_i\neq 0} \lambda_i$ endows $\Lambda_0$ with a valuation ring structure. Denote by $\Lambda_0^+$ the subset of $\Lambda_0$ consisting of elements with strictly positive valuation. This is the unique maximal ideal of $\Lambda_0$. In~\cite{FOOO} there was an additional parameter $e$ to encode Maslov index to have a $\Z$-graded $A_\infty$ structure. If we do not use this parameter, we need to work with $\Z/2\Z$-graded $A_\infty$ algebras. 

We briefly recall the construction of the Fukaya $A_\infty$ algebra structure on $\Omega^*(L,\Lambda_0)$. Let $\beta\in\pi_2(M,L)$ be a class in the relative homotopy group, and choose an almost complex structure $J$ compatible with (or simply tamed) $\omega$. Form $\cM_{k+1,\beta}(M,L;J)$, the moduli space of stable $(k+1)$-marked $J$-holomorphic disks in $M$ with boundary lying in $L_0$ of homotopy class $\beta$ with suitable regularity condition in interior and on the boundary. The moduli space $\cM_{k+1,\beta}(M,L;J)$ is of virtual dimension $d+k+\mu(\beta)-2$ (here $\mu(\beta)$ is the Maslov index of $\beta$).

There are $(k+1)$ evaluation maps $ev_i:\cM_{k+1,\beta}(M,L;J)\ra L$ for $i=0,\cdots,k$ which can be used to define a map $m_{k,\beta}:(\Omega^*(L,\C)^{\otimes k}) \ra \Omega^*(L,\C)$ of form degree $2-\mu(\beta)-k$  by formula
\[ m_{k,\beta} (\alpha_1,\cdots,\alpha_k):=(\ev_0)_!(\ev_1^*\alpha_1\wedge\cdots\wedge \ev_k^*\alpha_k).\]
To get an $A_\infty$ algebra structure we need to combine $m_{k,\beta}$ for different $\beta$'s. For this purpose we define a submonoid $G(L)$ of $\R^{\geq 0}\times 2\Z$ as the minimal one generated by the set 
\[\left\{(\int_\beta \omega,\mu(\beta))\in \R^{\geq 0}\times 2\Z \mid \beta\in \pi_2(M,L), \cM_{0,\beta}(M,L;J)\neq \emptyset\right\}.\]
Then we can define the structure maps $m_k:(\Omega(L,\Lambda_0))^{\otimes k} \ra \Omega(L,\Lambda_0)$ by
\[ m_k(\alpha_1,\cdots,\alpha_k):=\sum_{\beta\in G(L)} m_{k,\beta}(\alpha_1,\cdots,\alpha_k) T^{\int_\beta \omega}.\]
Note that we need to use the Novikov coefficients here since the above sum might not converge for a fixed value of $T$. The boundary stratas of $\cM_{k+1,\beta}(M,L;J)$ are certain fiber products of the diagram
\[\begin{CD}
\cM_{i+1,\beta_1}(M,L;J)\times_L\cM_{j+1,\beta_2}(M,L;J) @>>> \cM_{j+1,\beta_2}(M,L;J) \\
@VVV                      @VV\ev_l V\\
\cM_{i+1,\beta_1}(M,L;J) @>\ev_0>> L.
\end{CD}\]
Here $1\leq l\leq j$, $i+j=k+1$, and $\beta_1+\beta_2=\beta$. Indeed using this description of the boundary stratas the $A_\infty$ axiom for structure maps $m_k$ is an immediate consequence of Stokes formula.

We should emphasize that a mathematically rigorous realization of the above ideas involves lots of delicate constructions. Indeed the moduli spaces involved $\cM_{k+1,\beta}(M,L;J)$ are not smooth manifolds, but Kuranishi orbifolds with corners, which causes trouble to define an integration theory. Even if this regularity problem is taken of there are still transversality issues to define maps $m_{k,\beta}$ to have the expected dimension.  Moreover it is not enough to take care of each individual moduli space since the $A_\infty$ relations for $m_k$ follows from analyzing the boundary stratas in $\cM_{k+1,\beta}(M,L;J)$. Thus one needs to prove transversality of evaluation maps that are compatible for all $k$ and $\beta$. Furthermore one also need to deal with not only disk bubbles, but also sphere bubbles and regularity and transversality issues therein. We refer to the original constructions of~\cite{FOOO} and~\cite{Fukaya} for solutions of these problems. 

\paragraph{Geometric realization of pseudo-isotopies.} We shall apply the algebraic constructions in the previous section in the case $\overline{V}=\Omega^*(L,\C)$ is the de Rham differential graded algebra. In the construction of the Fukaya algebra of $L$ we need to make a choice of an $\omega$-tamed almost complex structure on $M$. One of the main results in~\cite{Fukaya} is that this $A_\infty$ algebra is uniquely determined up to pseudo-isotopies.

Indeed if $J_0$ and $J_1$ are two $\omega$-tamed almost complex structures we may choose a path $\cJ_t \; (t\in [0,1])$ of almost complex structures connecting them. This is always possible since the space of $\omega$-tamed almost complex structures is contractible, and hence path connected. Consider the parametrized moduli spaces
\[ \cM_{k+1,\beta} (M,L;\cJ):= \coprod_{t\in [0,1]} {t}\times \cM_{k+1,\beta}(M,L;\cJ_t).\]
Again we have $k$ evaluation maps
\[ \ev_i: \cM_{k+1,\beta}(M,L;\cJ) \ra L\]
for $1\leq i\leq k$. For the case $i=0$ we get another evaluation map
\[ \ev_0: \cM_{k+1,\beta}(M,L;\cJ)\ra L\times [0,1].\]
Then we can define an element 
\[\gamma:=\prod_{k,\beta} (m^t_{k,\beta}+h^t_{k,\beta})T^{\int_\beta \omega} \in \prod_k \Hom((V\otimes\Omega^*_{\Delta^1})^k,V\otimes\Omega^*_{\Delta^1}).\]
Here the structures maps $m^t_{k,\beta}$ and $h^t_{k,\beta}$ are defined by formula
\[ m^t_{k,\beta}(\alpha_1,\cdots,\alpha_k)+h^t_{k,\beta}(\alpha_1,\cdots,\alpha_k)dt:= (\ev_0)_![\ev_1^*\alpha_1\wedge\cdots\wedge\ev_k^*\alpha_k].\]
The following theorem of Fukaya verifies that $\gamma$ defined as above is indeed a pseudo-isotopy between the Fukaya algebras associated to $J_0$ and $J_1$. 
\begin{theorem}[~\cite{Fukaya} Section 11]
\label{thm:fukaya_isotopy}
The maps $m^t_{k,\beta}$ and $h^t_{k,\beta}$ defined above satisfies equation~\ref{eq:isotopy}.
\end{theorem}

\paragraph{Geometric realization of $2$-pseudo-isotopies.} In the same way if there is a $2$-simplex 
\[ \cJ_{x_0, x_1, x_2}: \Delta^2 \ra \cJ(M,\omega)\]
in the space $\cJ(M,\omega)$ of $\omega$-tamed almost complex structures bounding three paths of such. Here $x_0$, $x_1$, $x_2$ are coordinates on $\Delta^2$. Form the parametrized moduli spaces
\[ \cM_{k+1,\beta}(M,L;\cJ_{x_0, x_1, x_2}):=\coprod_{(x_0, x_1, x_2)\in \Delta^2} {(x_0, x_1, x_2)}\times \cM_{k+1,\beta}(M,L;\cJ_{(x_0, x_1, x_2)}).\]
Then using the evaluation maps 
\[\ev_i: \cM_{k+1,\beta}(M,L;\cJ_{x_0, x_1, x_2}) \ra L\]
for $1\leq i\leq k$ and 
\[\ev_0: \cM_{k+1,\beta}(M,L;\cJ_{x_0, x_1, x_2}) \ra L\times \Delta^2,\]
we can define a $2$-pseudo-isotopy whose component maps are given by
\[\alpha_1\otimes\cdots\otimes\alpha_k \mapsto (\ev_0)_![\ev_1^*\alpha_1\wedge\cdots\wedge\ev_k^*\alpha_k].\]
For purposes of this paper we only need to consider $n$-pseudo-isotopies for $n$ less or equal to $2$. In general existence of all the higher pseudo-isotopies are important, for instance to globalize the construction of the sheaf of symplectic functions in~\cite{Tu}. This aspect of pseudo-isotopies will be discussed in a forthcoming work~\cite{Tu2}.

\section{Reconstruction}
\label{sec:re}

This section contains the main result of this paper which describes a construction of a rigid analytic space $M^\chk_0$ associated to a Lagrangian torus fibration $\pi:M\ra B$ with certain assumptions. The space $M^\chk_0$ is constructed by gluing local pieces of affine rigid analytic spaces through transition functions which encodes instanton counting of pseudo-holomorphic disks. Moreover there is a canonical surjective map $M^\chk_0\ra B_0$ where $B_0$ denotes the smooth locus of $\pi$. Conjecturally this rigid analytic space $M^\chk_0$ should be an open dense part of the mirror manifold of $M$.

Throughout this section we shall assume that all smooth Lagrangian torus fibers of $\pi:M\ra B$ have Maslov index zero. We will also frequently use basic elements from the theory of rigid analytic spaces which we refer to the book~\cite{BGR} for a detailed discussion.

\paragraph{Open covering.} We need to work with an open covering of $B_0$ whose members are open subsets of $B_0$ additional properties. The construction of these open subsets is as follows.

Let $u\in B_0$ be a point in the smooth locus, and denote by $L_u$ the fiber Lagrangian $\pi^{-1}(u)$. Let $V\subset B_0$ be an contractible open neighborhood of $u$ such that there exists an action coordinates on $\pi^{-1}(V)$. We choose such a coordinate system $\phi: V\ra \R^n$ such that $\phi(u)=0$. Since $V$ is contractible the data $(u,V,\phi)$ induces an identification of symplectic manifolds
\[ s: \pi^{-1}(V) \ra L_u \times \phi(V) \subset T^*L_u=L_u\times H^1(L_u,\R)\]
through which $L_u\subset \pi^{-1}(V)$ is identified with $L_u\times 0$, i.e. the zero section of $T^*L_u$.

Let $r_1,\cdots, r_n\in \R$ be positive real numbers, and denote by $D_{r_1,\cdots,r_n}:= (-r_1,r_1)\times\cdots\times (-r_n,r_n)$ the corresponding poly-disk in $\R^n$. Let $U\subset V$ be an open subset of $V$ of the form $\phi^{-1}(D_{r_1,\cdots,r_n})$ for some $r_1,\cdots, r_n$.
For a point $p\in U$, let us define a \emph{diffeomorphism} $f_{u,p}$ of the symplectic manifold $M$ by formulas
\begin{align}
\phi(p) &=(p_1,\cdots, p_n)\in \R^n ,\\
V_{u,p} &:= p_1\partial/\partial x_1+\cdots+ p_n\partial/\partial x_n,\\
T_{u,p} &:= \epsilon \cdot s^* V_{u,p},\\
\label{eq:diff}
f_{u,p} &:= \mbox{time one flow of $T_{u,p}$}.
\end{align}
Here $\epsilon$ is a cut-off function supported in the neighborhood $V$, and is constant $1$ on $U$. The diffeomorphism $f_{u,p}$ has the property $f_{u,p}(L_u)=L_p$.

Finally by decreasing the positive numbers $r_1,\cdots, r_n$ we can assume that the diffeomorphism $f_{u,p}$ is tamed by the symplectic form $\omega$ for all $p\in U$. In the following we shall write the above data as a triple $(u,U, \phi)$ while not mentioning the set $V$ or the choice of $\epsilon$. We shall call such a triple $(u,U,\phi)$ an $\omega$-tamed triple.

\paragraph{Novikov field.} Let $\Lambda$ be the universal Novikov field with $\C$ coefficients. Explicitly this field is defined as
\[\Lambda:=\left\{\sum_{i=1}^\infty a_i T^{\lambda_i} \mid a_i\in \C, \lambda_i\in \R, \lim_{i\ra \infty} \lambda_i= \infty\right\}.\]
This is the field of fractions of the ring $\Lambda^\C_0$. It was proven in~\cite{FOOOtoric} that $\Lambda$ is algebraically closed. Of great importance for us is the existence of a valuation map $\val: \Lambda\ra \R$ defined by
\[ \val(\sum_{i=1}^\infty a_i T^{\lambda_i}):= \inf_i \lambda_i.\]

\paragraph{Local pieces.} Let $\pi:M\ra B$ be a Lagrangian torus fibration, and let $B_0\subset B$ be its smooth locus. Choose an open cover of $B_0$ by contractible open subsets $\left\{ U_i\right\}_{i\in\cI}$ where each $U_i$ is from a given $\omega$-tamed triple $(u_i, U_i, \phi_i)$. The fact that our valuation map $\val$ takes value in $\R$ ensures that all the open subsets $\phi_i(U_i)$ of $\R^n$ are rational domains. In fact all open subsets of the form $\phi_{i_j}(U_{i_1}\cap\cdots\cap U_{i_k})$ are rational domains for any $0\leq j\leq k$.

\begin{definition}
\label{def:star}
For a rational domain $U\subset \R^n$ we say a formal series of the form 
\[\sum_{(k_1,\cdots,k_n)\in \Z^n} a_{k_1\cdots k_n} (z^i_1)^{k_1}\cdots (z^i_n)^{k_n} \]
with $a_{k_1\cdots k_n}\in \Lambda$ satisfies \emph{condition $(*)$ of $U$} if the following holds.
\begin{equation}
\label{eq:*}
\val(a_{k_1\cdots k_n})+\sum_{j=1}^n k_j x_j \ra \infty \mbox{ \; as \;} (k_1,\cdots,k_n)\ra \infty, \forall (x_1,\cdots, x_n)\in U
\end{equation}
\end{definition}

Given a rational domain $U\subset \R^n$ the set of formal series that satisfies condition $(*)$ as in the above definition is a Tate algebra. The Tate algebra associated to an $\omega$-tamed triple $(u_i,U_i,\phi_i)$ shall be denoted by $\cO_i$. This is the set of formal series of the form 
\[\sum_{(k_1,\cdots,k_n)\in \Z^n} a_{k_1\cdots k_n} (z^i_1)^{k_1}\cdots (z^i_n)^{k_n}\]
such that $a_{k_1\cdots k_n}\in\Lambda$, and satisfies condition $(*)$ of $\phi_i(U_i)$. The fact that $\phi_i(U_i)$ is a rational domain ensures that $\cO_i$ is a Tate algebra. Taking the spectrum gives rise to a collection of rigid analytic spaces $\cU_i:=\Spec \cO_i\;(i\in \cI)$ over  the field $\Lambda$. These affine rigid analytic spaces will be our building blocks in the construction of $M^\chk_0$. For any pair of indices $i,\; j\in \cI$ the open subset $\phi_i(U_i\cap U_j)$ is also a rational domain. We write its associated Tate algebra by $\cO_{ij}$ and the corresponding rigid analytic space by $\cU_{ij}$.

\paragraph{Gluing data.} Recall a gluing data of the collection $\cU_i \;(i\in \cI)$ includes
\begin{itemize}
\item[(\rnum{1})] for any distinct pair of index $i$, $j\in \cI$, an open subset $\cU_{ij}\subset \cU_i$;
\item[(\rnum{2})] isomorphisms of analytic spaces $\Psi_{ij}: \cU_{ij}\ra \cU_{ji}$;
\item[(\rnum{3})] for each $i$, $j$, $\Psi_{ij}=\Psi_{ji}^{-1}$;
\item[(\rnum{4})] for each $i$, $j$, $k$, $\Psi_{ij}(\cU_{ij}\cap \cU_{ik})=\cU_{ji}\cap \cU_{jk}$, and $\Psi_{ik}=\Psi_{jk}\circ\Psi_{ij}$ on $\cU_{ij}\cap \cU_{ik}$.
\end{itemize}

\paragraph{Link to Maurer-Cartan spaces.} To define the transition functions we need to encode instanton corrections from symplectic geometry. For this consider the minimal model Fukaya algebra $H^*(L_{u_i},\Lambda)$. This vector space over $\Lambda$ is trivialized using the affine coordinates $\phi_i$ on $U_i$. Let us denote by $e^i_1,\cdots, e^i_n$ linear one-forms on $L_{u_i}$ associated to $\phi_i$ (the basis  for $H^*(L_{u_i})$ is then the wedge products of these one-forms). Observe there is an exponential map
\[ \exp: H^1(L_{u_i},\Lambda^+)\cong (\Lambda^+)^n \ra \Spec \cO_i\]
which in coordinates is given by $z^i_1=\exp( x^i_1),\cdots,z^i_n=\exp(x^i_n)$. Note that the exponential map is only well-defined on elements of positive valuation. 

The exponential map allows us to make connections with Maurer-Cartan moduli spaces associated to Fukaya algebras. Indeed there is a map $W: H^1(L_{u_i},\Lambda^+)\ra H^2(L_{u_i},\Lambda^+)$ defined by the Maurer-Cartan equation
\[ W: b\mapsto m_0+m_1(b)+\cdots+m_k(b^k)+\cdots.\]
By definition, the operator $m_k$ is a sum of the form $\sum_\beta m_{k,\beta}$. Moreover, by~\cite[Lemma 13.1]{Fukaya}, we have
\[ m_{k,\beta}(b^k)=\frac{1}{k!} m_{0,\beta},\]
which expresses certain compatibility of $m_k$ with the forgetful maps. Thus, for $b=x_1^i e^i_1+\cdots+x_n^i e^i_n$ we get
\begin{align*}
W(b) &=\sum_\beta\exp(\langle\partial \beta,b\rangle) m_{0,\beta} T^{\int_\beta \omega}\\
        &=\sum_\beta m_{0,\beta}(z^i_1)^{\langle\partial\beta,e^i_1\rangle}\cdots(z^i_n)^{\langle\partial\beta,e^i_n\rangle} T^{\int_\beta\omega}\\
        &=\sum_\beta m_{0,\beta} Z^i_\beta.
\end{align*}
Here in the last formula we introduced the notation 
\begin{equation}
\label{eq:zbeta}
Z^i_\beta:=(z^i_1)^{\langle\partial\beta,e^i_1\rangle}\cdots(z^i_n)^{\langle\partial\beta,e^i_n\rangle} T^{\int_\beta\omega}.
\end{equation}
It is important to observe that the above formula \emph{does not} immediately imply that the map $W:H^1(L_{u_i},\Lambda^+)\ra H^2(L_{u_i},\Lambda^+)$ factors through the exponential map $\exp: H^1(L_{u_i},\Lambda^+)\cong (\Lambda^+)^n \ra \Spec \cO_i$ because for that purpose we need to show that $W$ is actually convergent on entire $\Spec \cO_i$. This convergence problem was solved ingeniously by Fukaya in~\cite{Fukaya} which we now explain.

\paragraph{Fukaya's trick.} By definition of $\cO_i$, we need to show that the function $W$ is convergent for all points $(z^i_1,\cdots,z^i_n)\in (\Lambda^*)^n$ satisfying the property that $(\val(z^i_1),\cdots,\val(z^i_n))\in \phi_i(U_i)\subset \R^n$. For this purpose Fukaya invented a nice trick in~\cite{Fukaya}. The main idea is to identify the value of $W(z^i_1,\cdots,z^i_n)$ defined on the central fiber $L_{u_i}$ with a similar function defined on a near-by Lagrangian fibers. The convergence of the latter function will be automatic in view of Gromov's compactness theorem.

Let $(u_i,U_i,\phi_i)$ be an $\omega$-tamed triple. Let $p\in U_i$ be a point, and let $f_{u_i,p}$ be the diffeomorphism defined as in~\ref{eq:diff}. We consider the Fukaya algebra of the Lagrangian fiber $L_p$ whose structure maps are obtained by using the $\omega$-tamed almost complex structure $(f_{u_i,p})_*J$. The main advantage of this choice of almost complex structure is that the diffeomorphism $f_{u_i,p}$ induces an identification of various moduli spaces
\[ \cM_{k,\beta}(M,L_{u_i}; J) \stackrel{(f_{u_i,p})_*}{\longrightarrow} \cM_{k,(f_{u_i,p})_*\beta}(M,L_p; (f_{u_i,p})_*J)\]
that are compatible with all evaluation maps. This implies that the maps $m_{k,\beta}$ on $H^*(L_{u_i},\C)$ are the same as maps $m_{k,(f_{u_i,p})_*\beta}$ on $H^*(L_p,\C)$. Here note that the two vector spaces $H^*(L_{u_i},\C)$ and $H^*(L_p,\C)$ are both trivialized using affine coordinates $\phi_i$, and hence are isomorphic. Thus the structure maps $m^{u_i}_k:=\sum_{\beta}m_{k,\beta} T^{\int_\beta\omega}$ and $m^p_k:=\sum_{(f_{u_i,p})_*\beta} m_{k,(f_{u_i,p})_*\beta} T^{\int_{(f_{u_i,p})_*\beta}\omega}$ only differ through the part involves symplectic area. This difference can be made more precise through the following lemma.

\begin{lemma}
\label{lem:area}
Let the notations be as above. Then for all $\beta$ we have
\[ \int_{(f_{u_i,p})_*\beta} \omega-\int_\beta \omega = \langle\sum_{k=1}^n p^i_ke^i_k,\partial \beta\rangle\]
where $\partial\beta\in \pi_1(L_{u_i})$ is the boundary of $\beta$.
\end{lemma}

\medskip
\noindent This lemma implies that for all $\beta$ we have
\begin{equation}
\label{eq:maps}
m_{k,(f_{u_i,p})_*\beta}T^{\int_{(f_{u_i,p})_*\beta}\omega}=m_{k,\beta}T^{\int_\beta\omega}T^{\langle\sum_{k=1}^n p^i_ke^i_k,\partial \beta\rangle}.
\end{equation}
We are ready to prove the desired convergence property of $W$. Indeed let $(z^i_1,\cdots,z^i_n)\in (\Lambda^*)^n$ be such that the point $p:=(\val(z_1^i),\cdots,\val(z_n^i))\in \phi_i(U_i)$. 
The convergence of $W(z^i_1,\cdots,z^i_n)$ follows from the following computation using formula~\ref{eq:maps}.
\begin{align*}
 W_{u_i}(z^i_1,\cdots,z^i_n) &= \sum_\beta m_{0,\beta} Z_{\beta}\\
                &= \sum_{(f_{u_i,-p})_*\beta} m_{0,(f_{u_i,-p})_*\beta} Z_{(f_{u_i,-p})_*\beta}\cdot T^{\sum_{k=1}^n \langle-\val(z^i_k) e^i_k,\partial \beta\rangle}\\
               &=W_p (T^{-\val(z^i_1)}z^i_1,\cdots,T^{-\val(z^i_n)} z^i_n).
\end{align*}
Here the subscripts on $W$ is to indicate on which Lagrangian fiber we are applying the map $W$. In the last formula we are using the map $W$ on the Lagrangian torus fiber $L_p$ for $p:=(\val(z_1^i),\cdots,\val(z_n^i))\in \R^n$, which implies convergence since the point $T^{-\val(z^i_1)}z^i_1,\cdots,T^{-\val(z^i_n)} z^i_n$ lies in $(\Lambda_0)^n$. 

This proof of convergence in fact yields more. Let us denote by $\cO_i(p)$ the Tate algebra consisting of formal series $\sum_{(k_1,\cdots,k_n)\in \Z^n} a_{k_1\cdots k_n} (z^i_1)^{k_1}\cdots (z^i_n)^{k_n}$ which satisfies the condition $(*)$ for the rational domain $\phi_i(U_i)-\phi_i(p)$. Then we have the following commutative diagram 
\begin{equation}
\label{eq:potential}
\begin{CD}
\Spec \cO_i @> z^i_k\mapsto T^{-p^i_k} z^i_k>> \Spec \cO_i(p)\\
@V W VV      @V W VV\\
H^2(L_{u_i},\Lambda)  @>\cong >> H^2(L_p,\Lambda)
\end{CD}
\end{equation}
where the bottom horizontal map is the isomorphism induced by affine coordinates $\phi_i$.

\paragraph{Transition functions $\Psi_{ij}$.} In the previous paragraph we have introduced an isomorphism of rigid analytic spaces
\begin{equation}
\label{eq:translation}
S_{u_i,p}: \Spec \cO_i \ra \Spec\cO_i(p), \;\; z_k^i\mapsto T^{-p^i_k}z^i_k.
\end{equation}
If we take a point $p\in U_i\cap U_j$, then the image of $S_{u_i,p}$, when restricted to the analytic subspace $\cU_{ij}=\Spec \cO_{ij}\subset \cU_i=\Spec\cO_i$, is simply $\Spec\cO_{ij}(p)$ where the Tate algebra $\cO_{ij}(p)$ consists of formal power series satisfying condition $(*)$ of the domain $\phi_i(U_i\cap U_j)-\phi_i(p)$. To proceed further we need to make the following important assumption.

\begin{assumption}[Unobstructedness assumption]
\label{ass:unobs}
For any point $u_i$ in an $\omega$-tamed triple $(u_i,U_i,\phi_i)$, the map $W: H^1(L_{u_i},\Lambda^+)\ra H^2(L_{u_i},\Lambda^+)$ is identically zero.
\end{assumption} 

\noindent According to the diagram~\ref{eq:potential} this assumption implies the vanishing of $W_p$ for all $p\in U_i$. Equivalently this means that every element of $H^1(L_p,\Lambda^+)$ is a Maurer-Cartan element of $H^*(L_p,\Lambda^+)$. Moreover since the degree zero part of the Fukaya $A_\infty$ algebra $H^*(L_p,\Lambda^+)$ is strictly unital with $H^0(L_p,\Lambda)$ the subspace consisting of multiples of this unit, there is no non-trivial gauge equivalence defined on the set $H^1(L_p,\Lambda^+)$. Thus we conclude that
\[\MC^+(H^*(L_p,\Lambda^+))=H^1(L_p,\Lambda^+).\]
All the previous constructions also apply to the open subset $U_j$ endowed with coordinates $\phi_j$. However we have used two $A_\infty$ structures on $H^*(L_p,\Lambda^+)$ corresponding to the choice of two differential almost complex structures $(f_{u_i,p})_*J$ and $(f_{u_j,p})_*J$. As we have seen in the previous section, the two different $A_\infty$ structures are pseudo-isotopic, and a choice of a path in the space of $\omega$-tamed almost complex structures connecting the two points $(f_{u_i,p})_*J$ and $(f_{u_j,p})_*J$ specifies such an isotopy. There is in fact a natural choice. Recall by construction $f_{u_i,p}$ was defined as the time one flow of a vector field $T_{u_i,p}$. Denote this flow by $f_{u_i,p}^t$. Similarly denote by $f_{u_j,p}^t$ the flow of the vector field $T_{u_j,p}$. Thus the (smooth) concatenation $(f_{u_i,p}^{1-t})_*J\,\sharp\, (f_{u_j,p}^t)_*J$ of paths of almost complex structures satisfies the required boundary condition. By constructions in Sections~\ref{sec:algebra},~\ref{sec:gr} this defines a pseudo-isotopy $\gamma_{ij}$ (and $\gamma^\can_{ij}$) which produces an $A_\infty$ homomorphism
\[  F^\can_{ij}: H^*(L_p,\Lambda^+; (f_{u_i,p})_*J )\ra H^*(L_p,\Lambda^+; (f_{u_j,p})_*J ).\]
The map $ F^\can_{ij}$ further induces an isomorphism
\[ ( F^\can_{ij})_* : H^1(L_p,\Lambda^+)\ra H^1(L_p,\Lambda^+)\]
of the corresponding Maurer-Cartan moduli spaces. Our next goal is to solve the following commutative diagram.
\begin{equation}
\label{diag:Phi}
\begin{xy} 
(0,20)*+{H^1(L_p,\Lambda^+)}="a"; (30,20)*+{H^1(L_p,\Lambda^+)}="b";%
(0,0)*+{\Spec\cO_{ij}(p)}="c"; (30,0)*+{\Spec\cO_{ji}(p)}="d";
{\ar@{->}^\exp "a";"c"};{\ar@{->}^\exp "b";"d"}
{\ar@{->}^{( F^\can_{ij})_*} "a";"b"}
{\ar@{-->}^{\Phi_{ij}} "c";"d"};
\end{xy}
\end{equation}
Then using the maps $\Phi_{ij}$ the desired transition maps $\Psi_{ij}:\Spec\cO_{ij}\ra \Spec\cO_{ji}$ can be constructed easily by combining $\Phi_{ij}$ with the isomorphisms in formula~\ref{eq:translation}. The construction of $\Phi_{ij}$ will follow from the following two lemmas.

\begin{lemma}
\label{lem:exp}
The map $( F^\can_{ij})_*$ is of the form
\[ b\mapsto b+\sum_\alpha ( F^\can_{ij})_{0,\alpha} e^{\langle\partial\alpha,b\rangle} T^{\int_\alpha\omega}.\]
\end{lemma}

\proof To distinguish the usage of the de Rham model with the canonical model, we denote by $\gamma_{ij}^\dR$ and $\gamma_{ij}^\can$ the associated pseudo-isotopies on the corresponding models defined using the path $(f_{u_i,p}^{1-t})_*J\,\sharp\, (f_{u_j,p}^t)_*J$. Using continuous family version of multi-sections to perturb various moduli spaces as is done in~\cite{Fukaya} we can define the pseudo-isotopy $\gamma^\dR$ so that it satisfies
\[ \gamma^\dR_{k,\alpha}(b^k)=\frac{1}{k!}\langle\partial\alpha,b\rangle^k \gamma_{0,\alpha}\]
for any closed one-form $b$ on $L$ with positive valuation~\footnote{This result can also be extended to closed one-form with zero valuation by considering the usual topology on the coefficient field $\C$.}. This equation implies that $m^{\dR,t}_{k,\alpha}$ and $h^{\dR,t}_{k,\alpha}$, as components of $\gamma^{\dR}_{k,\alpha}$, both satisfy the same equation as above. Using summing over trees formulas~\ref{eq:can_1} and~\ref{eq:can_2} we conclude that the components of $\gamma^\can_{k,\alpha}$ which are $m^{\can,t}_{k,\alpha}$ and $h^{\can,t}_{k,\alpha}$ also satisfy this equation. To see this it suffices to observe that deleting a leave labeled by $b$ for a tree in $\cO(k+1,\alpha)$ produces trees in $\cO(k,\alpha)$, and moreover the $G$-decorations also matches since $\langle\partial\alpha,b\rangle$ is a linear function. Arguing in the same way using formula~\ref{eq:phi} we see that the $A_\infty$ homomorphism satisfies
\[ ( F^\can_{ij})_{k,\alpha}(b^k)= ( F^\can_{ij})_{0,\alpha} \frac{1}{k!}\langle\partial\alpha,b\rangle^k.\]
Summing over $k$ and $\alpha$ yields the formula in this Lemma. Finally we note that the linear term $b$ on the right hand side corresponds to the only tree with no interior vertices (the vertical line tree) which gives the identity map. The lemma is proved. \qed

\medskip
\noindent The above lemma implies certain exponential dependence of the map $( F^\can_{ij})_*$. Let us explicitly write down this map in coordinate $b=x^i_1 e^i_1+\cdots+x^i_n e^i_n$ we get
\[ x^i_k\mapsto x^i_k+\sum_{\alpha} \langle( F^\can_{ij})_{0,\alpha},e^i_k\rangle e^{\langle\partial\alpha,x^i_1 e^i_1+\cdots+x^i_n e^i_n\rangle} T^{\int_\alpha\omega}.\]
Using the exponential coordinates $z^i_k:=\exp(x^i_k)$ and notation~\ref{eq:zbeta} the above transformation becomes
\[ z^i_k \mapsto z^i_k \cdot e^{\sum_\alpha \langle( F^\can_{ij})_{0,\alpha}, e^i_k\rangle Z^i_\alpha}.\]
We need to show certain convergence property of the right hand side. This is accomplished in the following lemma whose proof is analogous to Fukaya's trick introduced earlier.
\begin{lemma}[Fukaya's trick on pseudo-isotopies]
\label{lem:convergence}
The map 
\[z^i_k \mapsto z^i_k \cdot e^{\sum_\alpha \langle( F^\can_{ij})_{0,\alpha}, e^i_k\rangle Z^i_\alpha}\]
as defined above is an automorphism of the Tate algebra $\cO_{ij}(p)$.

\end{lemma}

\proof The main point is to prove certain convergence property for the formal series $(\sum_\alpha \langle( F^\can_{ij})_{0,\alpha}, e^i_k\rangle Z^i_\alpha)$. We first show that for a point $(z^i_1,\cdots,z^i_n)\in (\Lambda^*)^n$ such that $(\val(z^i_1),\cdots,\val(z^i_n))\in \phi_i(U_i\cap U_j)-\phi_i(p)\subset\R^n$ this series is convergent. Denote by $q:=\phi_i^{-1}(\val(z^i_1)+p^i_1,\cdots,\val(z^i_n)+p^i_n)\in U_i\cap U_j$ the corresponding point in the intersection. 

Recall that the pseudo-isotopy $\gamma_{ij}$ which induces the $A_\infty$ homomorphism $ F^\can_{ij}$ was constructed using the path  $\cJ^t:=(f_{u_i,p}^{1-t})_*J\,\sharp\, (f_{u_j,p}^t)_*J$  of almost complex structures. The idea is to use a path of diffeomorphisms $\cF^t$ to push-forward the path $\cJ^t$ to obtain a path of almost complex structures connecting $(f_{u_i,q})_*J$ with $(f_{u_j,q})_*J$. Moreover  for each $t\in[0,1]$ we require that
\[ \cF^t(L_p)=L_q.\]
Thus $\cF^t$ maybe viewed as a parameter version of Fukaya's trick. The boundary condition of $\cF^t_*\cJ^t$ implies a boundary condition for $\cF^t$ which is simply
\[ \cF^0=f_{u_i,q}\circ f_{u_i,p}^{-1}, \mbox{\;\; and \;\;\;} \cF^1=f_{u_j,q}\circ f_{u_j,p}^{-1}.\]
We shall construct the path $\cF^t$ as the concatenation of three paths of diffeomorphisms $\cF_1$, $\cF_2$ and $\cF_3$ whose definitions are in order. To construct $\cF_1$ we define a one-parameter vector field $T^i_{tq-p}$ by 
\begin{align*}
 V^i_{tq-p} &:= (tq^i_1-p^i_1)\partial/\partial x_1^i+\cdots+(tq^i_n-p^i_n)\partial/\partial x_n^i\\
T^i_{tq-p} &:= \epsilon_i \cdot (s_i)^* V^i_{tq-p}
\end{align*}
in local affine coordinates. Denote by $f^i_{tq-p}$ the corresponding time one flow of $T^i_{tq-p}$. Then we define the path
\[ \cF^t_1:= f_{u_i,(1-t)q} \circ f^i_{tq-p}.\]
By definition we have $\cF^0_1=f_{u_i,q}\circ f_{u_i,p}^{-1}$ and $\cF^1_1=f^i_{q-p}$. The construction of $\cF_1$ is illustrated in the following picture.

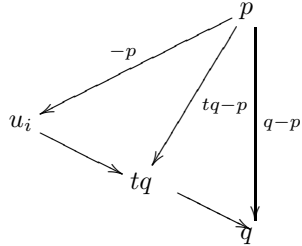
\[\begin{xy} 
(0,15)*+{u_i}="a"; (30,0)*+{q}="b";%
(30,30)*+{p}="c"; (16,7)*+{tq}="d";
{\ar@{->}_{-p} "c";"a"};{\ar@{->} "a";"d"}
{\ar@{->}^{tq-p} "c";"d"}
{\ar@{->}^{q-p} "c"; "b"}
{\ar@{->} "d"; "b"}
\end{xy}\]
Similarly working in the open subset $U_j$ we put
\begin{align*}
V^j_{(1-t)q-p} &:= ((1-t)q^j_1-p^j_1)\partial/\partial x^j_1 +\cdots+((1-t)q^j_n-p^j_n)\partial/\partial x^j_n\\
T^j_{(1-t)q-p} &:= \epsilon_{j}\cdot (s_j)^{*} V^j_{(1-t)q-p}.
\end{align*}
Denote by $f^j_{(1-t)q-p}$ the time one flow of $T^j_{(1-t)q-p}$. Define another path by
\[ \cF^t_3:= f_{u_j,tq}\circ f^j_{(1-t)q-p}.\]
We observe that $\cF^0_3=f^j_{q-p}$ and $\cF^1_3=f_{u_j,q}\circ f^{-1}_{u_j,p}$. This construction is illustrated in the picture below.

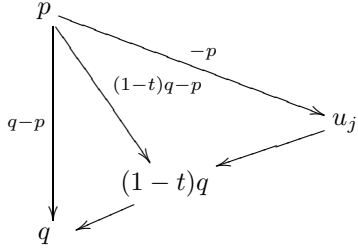
\[\begin{xy} 
(40,15)*+{u_j}="a"; (0,0)*+{q}="b";%
(0,30)*+{p}="c"; (16,7)*+{(1-t)q}="d";
{\ar@{->}^{-p} "c";"a"};{\ar@{->} "a";"d"}
{\ar@{->}^{(1-t)q-p} "c";"d"}
{\ar@{->}_{q-p} "c"; "b"}
{\ar@{->} "d"; "b"}
\end{xy}\]
Finally we need to construct a path connecting $f^i_{q-p}$ with $f^j_{q-p}$. For this we simply take the path
\[ \cF_2^t:=\mbox{\; time one flow of \;\;} (1-t)T^i_{q-p}+t T^j_{q-p}.\]
Then the smooth concatenation $\cF:=\cF_1\;\sharp\; \cF_2\;\sharp\; \cF_3$ satisfies all the desired conditions. Using $\cF$ we get a canonical identification of parametrized moduli spaces
\[ \cM_{k,\beta}(M,L_p; \cJ^t) \stackrel{(\cF^t)_*}{\longrightarrow} \cM_{k,(\cF^t)_*\beta}(M,L_q; (\cF^t)_*\cJ^t),\]
which implies that the associated isotopies $\gamma_{ij}$ and $\gamma_{ij}(q)$ only differ by the area term $T^{\int_\beta \omega}$ as in Lemma~\ref{lem:area}. This implies the sum
\[ \sum_\alpha \langle( F^\can_{ij})_{0,\alpha}, e^i_k\rangle Z^i_\alpha\]
when evaluated on the point $(z^i_1,\cdots,z^i_n)\in (\Lambda^*)^n$ with valuation vector given by $\phi_i(q)-\phi_i(p)=(q^i_1-p^i_1,\cdots,q^i_n-p^i_n)\in \R^n$ is equal to the function
\[ \sum_{(\cF^t)_*\alpha} \langle( F^\can_{ij}(q))_{0,(\cF^t)_*\alpha},e^i_k\rangle Z^i_{(\cF^t)_*\alpha}\]
evaluated on the point $(T^{p^i_1-q^i_1}z^i_1,\cdots, T^{p^i_n-q^i_n}z^i_n)$ whose valuation vector is $0$, hence convergent. This further implies that the term $\exp(\sum_\alpha \langle( F^\can_{ij})_{0,\alpha}, e^i_k\rangle Z^i_\alpha)$ remains to have valuation zero. Hence the map defined by
\[z_k^i\mapsto z_k^i\cdot e^{\sum_\alpha \langle( F^\can_{ij})_{0,\alpha}, e^i_k\rangle Z^i_\alpha}\]
is an automorphism of the Tate algebra $\cO_{ij}(p)$. The lemma is proved.\qed

\medskip
\remark The map $z_k^i\mapsto z_k^i\cdot e{\sum_\alpha \langle( F^\can_{ij})_{0,\alpha}, e^i_k\rangle Z^i_\alpha}$ has the following additional nice property. If $\beta\in H_2(M,L_p;\Z)$ is a relative homotopy class, then
\[ Z^i_\beta\mapsto Z^i_\beta \cdot e^{\sum_\alpha\langle( F^\can_{ij})_{0,\alpha},\beta\rangle Z^i_\alpha}.\]
This implies that this transformation is compatible with addition of homotopy classes. More precisely if we denote by $(Z_\beta^i)^*=Z^i_\beta \cdot e^{\sum_\alpha\langle( F^\can_{ij})_{0,\alpha},\beta\rangle Z^i_\alpha}$, then we have
\[ (Z^i_{\beta_1})^*\cdot(Z^i_{\beta_2})^*=(Z^i_{\beta_1+\beta_2})^*.\]
This form of transformation was suggested by Auroux in~\cite{Auroux}.

\paragraph{Explicit formula for transition functions.} To write down the map $\Phi_{ij}$ in the commutative diagram~\ref{diag:Phi} we need also to perform an affine change of ordinates since the right vertical exponential map in this diagram was written in coordinates $\phi_j$ rather than $\phi_i$. Thus if the two affine coordinates $x^i(p):=(x^i_1-p^i_1,\cdots,x^i_n-p^i_n)^t$, and $x^j(p):=(x^j_1-p^j_1,\cdots, x^j_n-p^j_n)^t$ transforms as
\[ x^j(p) = A\cdot x^i(p)\]
for some matrix $A=(A_{kl})\in \GL_n(\Z)$. In terms of exponential coordinates this gives
\[ z^j_k =\prod_l (z^i_l)^{A_{kl}}.\]
Hence the map $\Phi_{ij}$ can be written down in coordinates by
\begin{equation}
\label{eq:Phi}
z^j_k = \prod_l \;[z^i_l\cdot e^{\sum_\alpha \langle( F^\can_{ij})_{0,\alpha}, e^i_l\rangle Z^i_\alpha}]^{A_{kl}}.
\end{equation}
Lemma~\ref{lem:convergence} ensures the convergence property of $\Phi_{ij}$ (hence well-definedness), and Lemma~\ref{lem:exp} ensures the commutativity of diagram~\ref{diag:Phi}. Finally we define the transition functions $\Psi_{ij}: \Spec \cO_{ij} \ra \Spec \cO_{ji}$ by the composition
\begin{equation}
\label{eq:Psi}
 \Psi_{ij}:\Spec \cO_{ij} \stackrel{S_{u_i,p}}{\longrightarrow} \Spec \cO_{ij}(p) \stackrel{\Phi_{ij}}{\longrightarrow} \Spec\cO_{ji}(p) \stackrel{S_{u_j,p}^{-1}}{\longrightarrow} \Spec \cO_{ji}.
\end{equation}
Here the maps $S_{u_i,p}$ and $S_{u_j,p}$ are defined by formula~\ref{eq:translation}.

\paragraph{Gluing properties of transition functions.} Finally we prove that the transition maps defined above satisfy the desired gluing conditions.
\begin{lemma}
\label{lem:inverse}
We have $\Psi_{ij}\circ\Psi_{ji}=\id$.
\end{lemma}

\proof Since $\Psi_{ij}=S_{u_j,p}^{-1}\Phi_{ij}S_{u_i,p}$ we have
\[ \Psi_{ij}\circ\Psi_{ji}= S_{u_j,p}^{-1}\Phi_{ij}S_{u_i,p}S_{u_i,p}^{-1}\Phi_{ji}S_{u_j,p}=S_{u_j,p}^{-1}\Phi_{ij}\Phi_{ji}S_{u_j,p}.\]
Thus it suffices to show that $\Phi_{ij}\Phi_{ji}=\id$. By definition of $\Phi_{ij}$ there is a commutative diagram
\[\begin{CD}
H^1(L_p,\Lambda^+) @>(F^\can_{ji})_*>> H^1(L_p,\Lambda^+) @>( F^\can_{ij})_*>> H^1(L_p,\Lambda^+) \\
@V\exp VV @V\exp VV @V\exp VV\\
\Spec\cO_{ji}(p) @>\Phi_{ji}>> \Spec\cO_{ij}(p) @>\Phi_{ij}>> \Spec \cO_{ji}(p).
\end{CD}\]
The $A_\infty$ homomorphisms $ F^\can_{ij}$, $F^\can_{ji}$ are induced from inverse pseudo-isotopies, and hence their composition $ F^\can_{ij}\circ F^\can_{ji}$ is homotopic to the identity map. Taking the Maurer-Cartan functor yields a strict identification
\[ ( F^\can_{ij}\circ F^\can_{ji})_*=( F^\can_{ij})_*\circ(F^\can_{ji})_*=\id.\]
This implies the composition $\Phi_{ij}\circ\Phi_{ji}$ is also the identity map by formula in Lemma~\ref{lem:convergence}. \qed

To prove the cocycle condition for $\Psi_{ij}$'s we need the following lemma.

\begin{lemma}
\label{lem:ind}
The map $\Psi_{ij}$ constructed above is independent of the choice of $p\in U_i\cap U_j$.
\end{lemma}

\proof The proof is illustrated in the following commutative diagram.

\[\begin{xy} 
(0,20)*+{\Spec\cO_{ij}}="a"; (20,0)*+{\Spec\cO_{ij}(q)}="c"; (20, 40)*+{\Spec\cO_{ij}(p)}="b";%
(60,0)*+{\Spec\cO_{ji}(q)}="e"; (60,40)*+{\Spec\cO_{ji}(p)}="d"; (80, 20)*+{\Spec\cO_{ji}}="f";
{\ar@{->}^{S_{u_i,p}} "a";"b"};{\ar@{->}_{S_{u_i,q}} "a";"c"}; {\ar@{->}_{S^i_{p,q}} "b";"c"};
{\ar@{->}^{\Phi_{ij}(p)} "b";"d"};
{\ar@{->}^{\Phi_{ij}(q)} "c"; "e"};
{\ar@{->}^{S^j_{p,q}} "d"; "e"};
{\ar@{->}^{S_{u_j,p}^{-1}} "d"; "f"};
{\ar@{->}_{S_{u_j,q}^{-1}} "e"; "f"};
\end{xy}\]
The only non-trivial commutativity is the middle square, which follows from the Proof of Lemma~\ref{lem:convergence}. Indeed using the family of diffeomorphisms $\cF$ constructed there we can conclude that the two maps $\Phi_{ij}(p)$ and $\Phi_{ij}(q)$ only differ by the symplectic area term $\int_\alpha \omega$ which is corrected by the maps $S_{p,q}^i$ and $S_{p,q}^j$ by Lemma~\ref{lem:area}. The lemma is proved.\qed

\medskip
\noindent By the previous lemma we may pick a point $p\in U_i\cap U_j \cap U_k$ to prove the cocycle condition.
\begin{lemma}
\label{lem:image} Let the notations be as above. We have
\[ \Psi_{ij}(\Spec\cO_{ij} \cap \Spec\cO_{ik}) \subset \Spec\cO_{ji}\cap \Spec\cO_{jk}.\]
\end{lemma}

\proof To prove this we observe that a point $(z^i_1,\cdots,z^i_n)\in (\Lambda^*)n$ is inside $\Spec\cO_{ij} \cap \Spec\cO_{ik}$ if and only if its valuation vector $(\val(z^i_1,\cdots,z^i_n)$ lies inside $\phi_i(U_{ij}\cap U_{ik})$. Then we have
\begin{align*}
&\val(S_{u_i,p}(z^i_1,\cdots,z^i_n))= (\val(z^i_1)-p^i_1,\cdots,\val(z^i_n)-p^i_n);\\
&\val(\Phi_{ij}S_{u_i,p}(z^i_1,\cdots,z^i_n))=\\
&= (\sum_k A_{1k} (\val(z^i_k)-p^i_k),\cdots,\sum_k A_{nk} (\val(z^i_k)-p^i_k));\\
&\val(S^{-1}_{u_j,p}\Phi_{ij} S_{u_i,p}(z^i_1,\cdots,z^i_n))=\\
&=(p^j_1+\sum_k A_{1k} (\val(z^i_k)-p^i_k),\cdots,p^j_n+\sum_k A_{nk} (\val(z^i_k)-p^i_k)).
\end{align*}
Here the matrix $A\in \GL_n(\Z)$ is linear part of change of coordinates from $x^i$ to $x^j$. From this computation we see that if $(\val(z^i_1,\cdots,z^i_n)\in\phi_i(U_{ij}\cap U_{ik})$ then $\val(\Psi_{ij}(z^i_1,\cdots,z^i_n))\in \phi_j(U_{jk}\cap U_{ji})$. Thus we have proved that the lemma.\qed

\begin{lemma}
\label{lem:cocycle}
The cocycle condition $\Psi_{jk}\Psi_{ij}=\Psi_{ik}$ holds.
\end{lemma}

\proof By the construction of $\Psi$ we have
\begin{align*}
\Psi_{jk}\circ\Psi_{ij}&= S_{u_k,p}^{-1}\Phi_{jk}S_{u_j,p}S_{u_j,p}^{-1}\Phi_{ij}S_{u_i,p}\\
                    &= S_{u_k,p}^{-1}\Phi_{jk}\Phi_{ij} S_{u_i,p}.
\end{align*}
Hence it suffice to prove that $\Phi_{jk}\Phi_{ij}=\Phi_{ik}$. To this end let us denote by $\cU_{ij}(p):=\Spec\cO_{ij}(p)$ (similarly for other indices), and consider the following commutative diagram.

\[\begin{xy}
(0,0)*+{\cU_{ij}(p)\cap\cU_{ik}(p)}="a";
(40,0)*+{\cU_{ji}(p)\cap\cU_{jk}(p)}="b"; 
(80,0)*+{\cU_{ki}(p)\cap\cU_{kj}(p)}="c";
(0,30)*+{H^1(L_p,\Lambda^+)}="d";
(40,30)*+{H^1(L_p,\Lambda^+)}="e";
(80,30)*+{H^1(L_p,\Lambda^+)}="f";
{\ar@{->}^{( F^\can_{ij})_*} "d"; "e"};
{\ar@{->}^{(F^\can_{jk})_*} "e"; "f"};
{\ar@{->}^{\Phi_{ij}} "a"; "b"};
{\ar@{->}^{\Phi_{jk}} "b"; "c"};
{\ar@{->}^{\exp} "d"; "a"};
{\ar@{->}^{\exp} "e"; "b"};
{\ar@{->}^{\exp} "f"; "c"};
\end{xy}
\]
We claim that the top row composition of the above diagram is  \[(F^\can_{jk})_*\circ( F^\can_{ij})_*=(F^\can_{ik})_*.\]
To see this we use the fact that the space of $\omega$-tamed almost complex structures is contractible. Hence there exist a two-parameter family of almost complex structures $\cJ_{x_0, x_1, x_2} \; ((x_0,x_1,x_2)\in \Delta^2)$ such that
\begin{align*}
\cJ_{0,x_1,x_2}&=\cJ_{jk};\\
\cJ_{x_0,0,x_2}&=\cJ_{ik};\\
\cJ_{x_0,x_1,0}&=\cJ_{ij}.
\end{align*}
Here $\cJ_{ij}$ ($\cJ_{ik}$ and $\cJ_{jk}$ respectively) is the path of almost complex structures used to construct the homotopy $ F^\can_{ij}$ ($F^\can_{ik}$ and $F^\can_{jk}$ respectively). This two parameter family $\cJ_{x_0,x_1,x_2}$ induces a $2$-homotopy that bounds $1$-homotopies $\gamma_{ij}$, $\gamma_{jk}$ and $\gamma_{ik}$. Thus by Theorem~\ref{thm:homotopy} and homotopy invariance of Maurer-Cartan functor we get
$(F^\can_{jk})_*\circ( F^\can_{ij})_*=(F^\can_{ik})_*$. To this end we observe there is another commutative diagram
\[\begin{CD}
H^1(L_p,\Lambda^+) @>(F^\can_{ik})_*>> H^1(L_p,\Lambda^+)\\
@V\exp VV     @V\exp VV\\
\cU_{ij}(p)\cap\cU_{ik}(p) @>\Phi_{ik}>> \cU_{ki}(p)\cap\cU_{kj}(p).
\end{CD}\]
Combined with the previous diagram, it implies the cocycle condition
\[\Phi_{jk}\Phi_{ij}=\Phi_{ik}\]
for $\Phi$, which implies that of $\Psi$. The proof is complete.\qed

\medskip
\noindent It follows from Lemma~\ref{lem:inverse}, Lemma~\ref{lem:image} and Lemma~\ref{lem:cocycle} the local affinoid domains $\cU_i$ can be glued to obtain a rigid analytic space which we shall denote by $M^\chk_0$. Moreover the natural valuation maps $\val:\cU_i\ra \phi_i(U_i)$ are compatible with the gluing maps $\Psi_{ij}$, which follows from the computation in the proof of Lemma~\ref{lem:image}. Hence we get a global valuation map
\[ \val : M^\chk_0 \ra B_0 .\]
This might be thought of as the Non-Archimedean version of a torus fibration. We summarize the main results of this paper in the following corollary.
\begin{corollary}
\label{coro:main}
Let $\pi: M\ra B$ be a Lagrangian torus fibration with vanishing Maslov index. Assume furthermore Assumption~\ref{ass:unobs} holds. Then the local pieces $\cU_i:=\Spec\cO_i$ of affinoid domains can be glued via transition maps $\Psi_{ij}$ (see formula~\ref{eq:Psi}) to form a rigid analytic space $M^\chk_0$ over the universal Novikov field $\Lambda$. Moreover there is a globally defined valuation map
\[ \val: M^\chk_0\ra B_0.\]
\end{corollary}

\end{document}